\documentclass[twoside,11pt]{amsart} 
\usepackage{epsfig}
\usepackage{graphicx}
\usepackage{latexsym} 
\usepackage{amsmath} \usepackage{amssymb}
\usepackage{cancel}



\newtheorem{theorem}{Theorem}[section]
\newtheorem{lemma}[theorem]{Lemma}

\newdimen\epsfxsize
\newdimen\epsfysize

\newcommand {\gap}     {\makebox[0.075 in]{}}

\newcommand {\fto}     {\longrightarrow}

\newcommand {\set}[1]  {\left\{ {#1} \right\}}   
   
\newcommand {\cons}[1] {\langle {#1} \rangle }

\newcommand{\pml}[1] {{#1}^\#}

\newcommand{\Z}     {{\mathbb Z}}

\newcommand{\pgap}   {{\mathcal G}}

\newcommand{\wgp}[2] {\left. \bar{w}_{{#1}} \right|_{{#2}}}

\newcommand{\Bi}[2]{\left( \begin{array}{c}{#1} \\ {#2} \end{array}\right)}
\newcommand{\lil}{\scriptstyle}

\begin{document}

\title{On the last digits of consecutive primes}

\date{9 July 2016}

\author{Fred B. Holt}
\address{fbholt62@gmail.com ; 5520 - 31st Ave NE, Seattle, WA 98105}

\begin{abstract}
Recently Oliver and Soundararajan made conjectures based on computational enumerations about the 
frequency of occurrence of pairs of last digits for consecutive primes.
By studying Eratosthenes sieve, we have identified discrete dynamic systems that
exactly model the populations of gaps across stages of Eratosthenes sieve.
Our models provide some insight into the observed biases in the occurrences of last digits in consecutive primes,
and the models suggest that the biases will ultimately be reversed for large enough primes.

The exact model for populations of gaps across stages of Eratosthenes sieve provides
a constructive complement to the probabilistic models rooted in the work of Hardy and Littlewood.
\end{abstract}

\maketitle

\section{Introduction}

Recently Oliver and Soundararajan \cite{OS, OSQ} computed the 
distribution of the last digits of consecutive primes for the first $10^8$
prime numbers.  Their calculations revealed a bias:  the pairs $(1,1)$, $(3,3)$, $(7,7)$
and $(9,9)$ occur about a third less often than other ordered pairs of last digits 
of consecutive primes.  Their calculations are shown
in Table~\ref{BiasTable}.

For the past several years we have been studying the cycle of gaps
$\pgap(\pml{p})$ that
arises at each stage of Eratosthenes sieve.  Our work to this point is
summarized in \cite{HRSFU}.  

We have identified a population model that describes the growth of 
the populations of any gap $g$ in the cycle of gaps $\pgap(\pml{p})$, across
the stages of Eratosthenes sieve.
The recursion from one cycle of gaps,
$\pgap(\pml{p_{k-1}})$, to the next, $\pgap(\pml{p_k})$, leads to a discrete dynamic model
that provides exact populations for a gap $g$ in the cycle $\pgap(\pml{p})$, provided
that $g < 2p$.  The model provides precise asymptotics for the ratio of the population of the gap $g$ to
the population of the gap $2$ once the prime $p$ is larger than any prime factor of $g$.
This discrete dynamic system is deterministic, not probabilistic.  

The discrete dynamic model provides some insight into the phenomenon that
Oliver and Soundararajan have observed \cite{OS, OSQ}. 
\begin{enumerate}
\item We look at the asymptotic ratios of the populations of small gaps to the gap $g=2$.
These asymptotic ratios suggest that the reported biases will erode away for samples of
much larger primes.
\item We look at additional terms in the model, to understand rates of convergence to the
asymptotic values.
To first order this explains some of the biases exhibited in Table~\ref{BiasTable}.
\item We initially work in base $10$, so we then examine the results for 
a few different bases, to see how the biases depend on the base.
\end{enumerate}

\begin{figure}[t]
\centering
\includegraphics[width=5in]{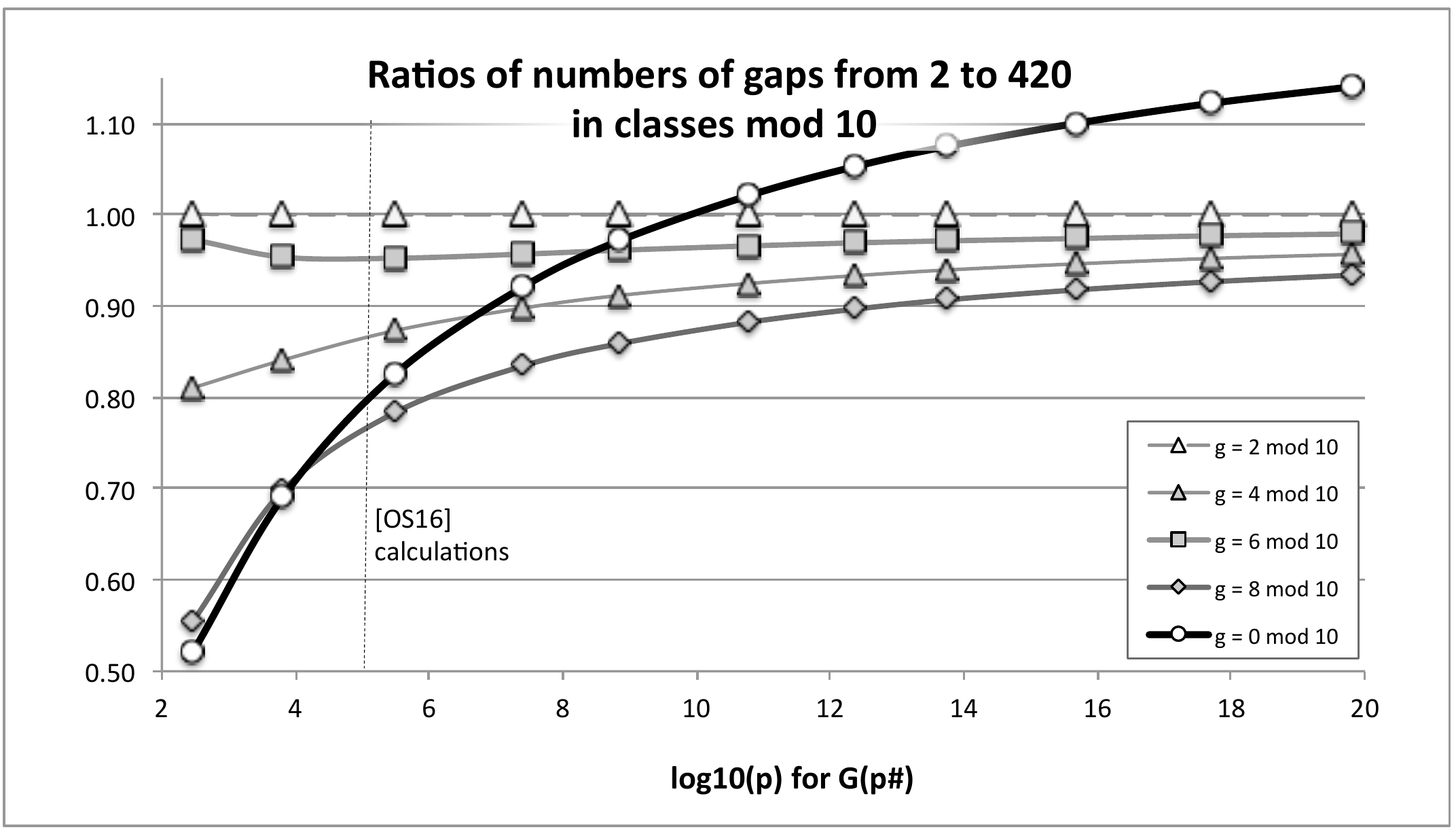}
\caption{\label{WratFig} A graph of the ratios of the populations of gaps in each residue class modulo $10$,
normalized by the population of gaps $g = 2 \bmod 10$.
Here the ratios in $\pgap(\pml{p})$ are approximated by Equation~\ref{EqEigSys} to twelve terms.
We used initial conditions from $\pgap(\pml{37})$ for gaps up to $g \le 420$.  
The dashed line indicates where
the calculations by Oliver and Soundararajan lie. }
\end{figure}

These observations apply to the stages of Eratosthenes sieve as the sieve proceeds.
All gaps between prime numbers arise in a cycle of gaps.  To connect our results 
to the desired results on gaps between primes, we would need to better understand how
gaps survive later stages of the sieve, to be affirmed as gaps between
primes.
Until the models for survival have a higher accuracy,
the results based on the exact models for $\pgap(\pml{p})$
can only be approximately applied to gaps between prime numbers.

We offer the exact model on populations of gaps in $\pgap(\pml{p})$ as a constructive
complement to the approaches working from the probabilistic models
pioneered by Hardy and Littlewood \cite{HL}.

\begin{table}
\begin{center}
\begin{tabular}{|ccr|ccr|} \hline
$a$ & $b$ & $\pi(x_0; 10, (a,b))$ & $a$ & $b$ & $\pi(x_0; 10, (a,b))$ \\ \hline
 $1$ & $1$ & $4,623,042$ & $7$ & $1$ & $6,373,981$ \\
  & $3$ & $7,429,438$ & & $3$ & $6,755,195$ \\
  & $7$ & $7,504,612$ & & $7$ & $4,439,355$ \\
  & $9$ & $5,442,345$ & & $9$ & $7,431,870$ \\ 
 $3$ & $1$ & $6,010,982$ & $9$ & $1$ & $7,991,431$ \\
  & $3$ & $4,442,562$ & & $3$ & $6,372,941$ \\
  & $7$ & $7,043,695$ & & $7$ & $6,012,739$ \\
  & $9$ & $7,502,896$ & & $9$ & $4,622,916$  \\  \hline
\end{tabular}
\caption{\label{BiasTable} Oliver and Soundararajan's table \cite{OS} of computed distributions
of last digits of consecutive primes for the first $10^8$ primes.  Here they are working
in base $10$.  In Section~\ref{SecBases} we address their calculations in base $3$ as well.}
\end{center}
\end{table}
 
 \subsection{Notes on version 3.}
 Since posting the first version of this work, we extended our analysis of the polynomial approximations in 
 section~\ref{AsymSection}.  We had initially worked with six terms of the polynomial expansion, and we 
have extended this to twelve terms.  This has helped us refine our claims about the convergence for 
the populations of some of the gaps.  Two examples of the progressive approximations, for the gaps $g=30$
and $g=420$, are shown in Figure~\ref{LPolyFig}.  

We have also introduced Mertens' Third Theorem, which ties the system parameter $\lambda = a_2^k$
in Equation~\ref{EqEigSys} to the magnitude of the prime $p$ in the cycle $\pgap(\pml{p})$.
 
 
 \section{The model for populations of gaps in $\pgap(\pml{p})$}
 Here we restate only a select few results from \cite{HRSFU} that are relevant to studying the
 last digits of consecutive primes.

We do not study the gaps between primes directly.  Instead, we study the cycle of gaps
$\pgap(\pml{p})$ at each stage of Eratosthenes sieve.  Here, $\pml{p}$ is the
{\em primorial} of $p$, which is the product of all primes from $2$ up to and including $p$.
$\pgap(\pml{p})$ is the cycle of gaps among the generators of $\Z \bmod \pml{p}$.
These generators and their images through the counting numbers are the candidate primes
after Eratosthenes sieve has run through the stages from $2$ to $p$.  All of the remaining primes
are among these candidates.

The cycle of gaps $\pgap(\pml{p})$
consists of $\phi(\pml{p})$ gaps that sum to $\pml{p}$.  For example, we have
$\pgap(\pml{5}) = 64242462$, and
$$\pgap(\pml{7}) =
 {\it 10}, 2424626424 {\it 6} 62642 {\it 6} 46 {\it 8} 42424 {\it 8} 64 {\it 6} 24626 {\it 6} 
  4246264242, {\it 10}, 2
$$

There is a substantial amount of structure preserved in the cycle of gaps from
one stage of Eratosthenes sieve to the next, from $\pgap(\pml{p_k})$ to 
$\pgap(\pml{p_{k+1}})$.  This structure is sufficient to enable us to give exact counts for gaps and for 
sufficiently short constellations in $\pgap(\pml{p})$ across all stages of the sieve.
 
As the sieve proceeds from one prime $p_{k-1}$ to the next $p_k$, the recursive construction 
leads to a discrete dynamic system that provides exact counts
of a gap and its driving terms.  The driving terms for a gap $g$ are constellations $s$ in $\pgap(\pml{p})$ that
have the same sum $g$.  For example, the driving terms for $g=6$ are the gaps $g=6$ themselves and the
constellations of length $2$ that sum to $6$:  e.g., $s=24$ and $s=42$.  
Under the closures in the recursion --  which correspond to eliminating a candidate from later stages of the sieve --
these constellations
will produce more gaps $g=6$.

These raw counts for populations of gaps grow
superexponentially by factors of $(p_k-2)$, and so to better understand their behavior we take the ratio
of a raw count $n_{g,j}(\pml{p})$ of the driving terms of length $j$ for the gap $g$ in $\pgap(\pml{p})$ 
to the number of gaps $g=2$ at this stage of the sieve.
\[ w_{g,j}(\pml{p_k}) = \frac{n_{g,j}(\pml{p_k})}{n_{2,1}(\pml{p_k})} \]

For a gap $g$ that has driving terms of lengths up to $j$, we take any $J \ge j$, and we form a vector
of initial values $\left. \bar{w}_g\right|_{\pml{p_{0}}}$, whose $i^{\rm th}$ entry is the ratio
$w_{g,i}$.  These values evolve according to the dynamic system:
\begin{eqnarray*}
\bar{w}_g (\pml{p_k}) & = & \left. M_{1:J} \right|_{p_k} \cdot \bar{w}_g(\pml{p_{k-1}}) \\
   & = & \left. M_{1:J} \right|_{p_k} \left. M_{1:J} \right|_{p_{k-1}} \cdots \left. M_{1:J} \right|_{p_1} \cdot \bar{w}_g(\pml{p_{0}}) \\
   & = &  M_{1:J}^k \cdot \bar{w}_g(\pml{p_{0}})
\end{eqnarray*}
Here $M(p)$ is a banded matrix with diagonal entries $a_j(p) = (p-j-1)/(p-2)$ and superdiagonal entries 
$b_j(p) = j / (p-2)$.  

We use the notation $M^k$ and $a_j^k$ to indicate the products over the range of primes
from $p_1$ to $p_k$, relative to some starting value $p_0$.
\begin{eqnarray*}
 M^k & = & \left. M \right|_{p_k} \cdot \left. M \right|_{p_{k-1}} \cdots \left. M \right|_{p_1} \\
 a_j^k & = & \prod_{p = p_1}^{p_k} \frac{p-j-1}{p-2}
\end{eqnarray*}
 
\subsection{Eigenstructure of the system matrix $M_{1:J}$}
The system matrix $M_{1:J}(p)$ is diagonalizable with a particularly nice eigenstructure.  
\[
M_{1:J}(p)  =  R \cdot \Lambda \cdot L
\]
with $LR = I$.

The right eigenvectors are the columns of $R$, and $R$ is an upper triangular Pascal matrix of
alternating sign:
\begin{eqnarray*}
R_{ij} & = & \left\{ \begin{array}{cl}
(-1)^{i+j}\Bi{j-1}{i-1} & {\rm if} \; i \le j \\ & \\
0 & {\rm if} \; i > j \end{array} \right.
\end{eqnarray*}

The left eigenvectors are the rows of $L$, and $L$ is an upper triangular Pascal matrix:
\begin{eqnarray*}
L_{ij}  & = & \left\{ \begin{array}{ll}
\Bi{j-1}{i-1} & {\rm if} \; i \le j \\ & \\
0 & {\rm if} \; i > j \end{array} \right. 
\end{eqnarray*}

The eigenvalues are the system coefficients $a_j = (p-j-1) / (p-2)$:
\[ \Lambda =  {\rm diag}(1, a_2, \ldots, a_J). \]

While the eigenvalues depend on
the prime $p$, the eigenvectors do not.  We thereby get a similarly nice eigenstructure for $M^k_{1:J}$.
\[
M^k_{1:J}   =  R \cdot \Lambda^k \cdot L
\]
in which the eigenvalues $\lambda^k_j = a_j^k = \prod_{p=p_1}^{p_k} (p-j-1) / (p-2)$.

\subsection{Implications of the discrete dynamic system}
We borrow a few results from \cite{HRSFU} that have direct bearing on the distributions
of last digits between consecutive primes.

Although the asymptotic growth of all gaps
is equal, the initial conditions and driving terms are important.
Brent \cite{Brent} made analogous observations.
His Table 2 indicates the importance of the lower-order effects in
estimating relative occurrences of certain gaps.

\begin{lemma}\label{LemmaEig}
For a gap $g$, let $p_0$ be any prime greater than the greatest prime factor of $g$,
and let $J$ be at least as large as the longest driving term for $g$ in $\pgap(\pml{p_0})$.
Then in $\pgap(\pml{p_k})$,
\begin{eqnarray}\label{EqEigSys}
w_{g,1}({\pml{p_k}}) & = &  (L_1 \wgp{g}{\pml{p_0}}) - a_2^k (L_2 \wgp{g}{\pml{p_0}}) \nonumber \\
 & & \gap + a_3^k (L_3 \wgp{g}{\pml{p_0}}) \cdots 
  + (-1)^{J+1} a_J^k (L_J \wgp{g}{\pml{p_0}}) \nonumber \\
  & \approx & (L_1 \wgp{g}{\pml{p_0}}) - a_2^k (L_2 \wgp{g}{\pml{p_0}})  \\
  & & \gap  + (a_2^k)^2 (L_3 \wgp{g}{\pml{p_0}}) \cdots 
  + (-1)^{J+1} (a_2^k)^{J-1} (L_J \wgp{g}{\pml{p_0}}) \nonumber
\end{eqnarray}
\end{lemma}

This is Equations~(6 \& 7) in \cite{HRSFU}.  From this expansion we can compute the asymptotic
values of the ratios $w_{g,1}(\pml{p})$, and we can analyze the rate of convergence to the asymptotic value.
To obtain the asymptotic ratio $w_{g,1}(\infty)$, since $L_1 = [\cons{1}]$ we simply add together the initial 
ratios of all driving terms.  As quickly as $a_2^k \fto 0$, the ratios $w_{g,1}(\pml{p_k})$ converge to the
asymptotic ratio $w_{g,1}(\infty)$.

For these asymptotic ratios, we restate Corollary~5.4 and Theorem~5.5 of \cite{HRSFU} here.

\begin{theorem}\label{PolThm}
For any $g=2n$, the gap $g$ eventually occurs in Eratosthenes sieve.
Let $\bar{q}$ be the largest prime factor of $g$.  Then for $p \ge \bar{q}$,
\[ w_{g,1}(\infty) = L_1 \wgp{g}{\pml{p}} = \sum w_{g,j}(\pml{p}) = \prod_{q > 2, \; q | g} \frac{q-1}{q-2} \]
\end{theorem}

This theorem establishes an analogue of Polignac's conjecture for 
Eratosthenes sieve \cite{HRPol}, that for any number $2n$ the gap $g=2n$ does occur infinitely often
in the sieve, and further that the ratio of occurrences of this gap to the gap $2$ approaches
the ratio implied by Hardy \& Littlewood's Conjecture B \cite{HL}.

\subsection{Estimating the rate of convergence of $a_k^2 \fto 0$}.  We have calculated $a_2^k$
for primes into the range of $10^{15}$, at which $a_k^2 \approx 0.105$ (with $p_0=37$).  
Even into this range, gaps of sizes in the low hundreds are still appearing in ratios far below their
asymptotic values.  We need a way to estimate $a_2^k$ for much larger primes.  

Mertens' Third Theorem provides these estimates.  The theorem is that
\[ \prod_{p \le q} \left( \frac{p-1}{p} \right) = \frac{e^{-\gamma}+o(1)}{\ln q}. \]
We define the constant $c_0 = \prod_{q \le p_0} q/(q-1)$.  Then we can establish
both a lower and an upper bound for $a_2^k$ for large primes $p_k$:
\begin{equation}\label{EqBnds}
\frac{p_0-1}{p_0} \cdot c_0 \cdot \frac{e^{-\gamma} + o(1)}{\ln p_{k-1}}  \; < \; a_2^k
 \; < \; c_0 \cdot \frac{e^{-\gamma} + o(1)}{\ln p_k} 
\end{equation}
For the calculations in this paper, we are using $p_0=37$.  We used these bounds to extend
Figure~\ref{WratFig} across the range $p_k \in [10^{15}, 10^{20}]$.


\section{Ultimate distributions of last digits of consecutive primes}
Consider an ordered pair $(a,b)$ of last digits of consecutive primes \cite{OS}, with $a,b \in \set{1,3,7,9}$.  
We are interested in the size of the
set of indices $k$, such that $p_k = a \bmod 10$ and $p_{k+1} = b \bmod 10$.

By way of example, suppose $b-a = 0 \bmod 10$.  Then $p_{k+1}-p_k = 10$, or
$p_{k+1}-p_k = 20$, or in general $p_{k+1}-p_k$ is some multiple of $10$.  So how often do these
gaps $g = 10, 20, 30, \ldots$ arise?

The ordered pairs $(a,b)$ of last digits correspond to specific gaps as follows
\begin{center}
\begin{tabular}{lcl}
\multicolumn{1}{c}{$(a,b)$'s} & $\Leftrightarrow $ & \multicolumn{1}{c}{$g$'s} \\ \hline
$(1,1), \; (3,3), \; (7,7), \; (9,9)$ & & $10, 20, 30, 40, \ldots $ \\
$(1,3), \; (7,9), \; (9,1)$ & & $2, 12, 22, 32, 42, \ldots $ \\
$(3,7), \; (7,1), \; (9,3)$ & & $4, 14, 24, 34, 44, \ldots $ \\
$(1,7), \; (3,9), \; (7,3)$ & & $6, 16, 26, 36, 46, \ldots $ \\
$(1,9), \; (3,1), \; (9,7)$ & & $8, 18, 28, 38, 48, \ldots $
\end{tabular}
\end{center}

This table already provides us with a couple of insights into the problem.
The class $b-a = 0 \bmod 10$ has four
ordered pairs and the other classes have three.  So in order for these ordered pairs to occur
equally often, the gaps in $g = 0 \bmod 10$ must occur $4/3$ as often as the gaps in the other
classes.  We also note that within any of the five classes, if the distribution of a single gap 
across the corresponding ordered pairs is not uniform, then this would lead to a biased distribution
within this class.

\renewcommand\arraystretch{0.8}
\begin{table}
\begin{center}
\begin{tabular}{|r|rrrr|rr|} \hline
gap & \multicolumn{4}{c|}{ $n_{g,j}(\pml{37})$: driving terms of length $j$ in $\pgap(\pml{37})$ } & 
$w_{g,1}(\pml{37})$ & $w_{g,1}(\infty)$ \\ [1 ex]
$g$ & $j=1$ & $2$ & $3$ & $4$ &  & \\ \hline 
 $\lil 2, \; 4$ & $\lil 217929355875$ & & & & $\lil{1}$ & $\lil{1}$ \\
 $\lil 6$ & $\lil 293920842950$ & $\lil 141937868800$ & & & $\lil 1.348698$ & $\lil{2}$ \\
 $\lil 8$ & $\lil 91589444450$ & $\lil 110741954050$ & $\lil 15597957375$ & & $\lil 0.420271$ &  $\lil{1}$ \\
 $\lil 10$ & $\lil 108861586050$ & $\lil 150514973700$ & $\lil 31195914750$ & & $\lil 0.499527$ & $\lil{4/3}$ \\
 $\lil 12$ & $\lil 83462164156$ & $\lil 219604134932$ & $\lil 121198832118$ & $\lil 11593580544$  & $\lil 0.382978$ & $\lil{2}$ \\
 $\lil 14$ & $\lil 83462164156$ & $\lil 115853913448$ & $\lil 93409823052$ & $\lil 17390370816$  & $\lil 0.159965$ & $\lil{6/5}$ \\
 $\lil 16$ & $\lil 16996070868$ & $\lil 78769359396$ & $\lil 91933104354$ & $\lil 28714181132$ & $\lil 0.077989$ &  $\lil{1}$ \\
 $\lil 18$ & $\lil 21218333416$ & $\lil 122467715552$ & $\lil 191942799048$ & $\lil 91130022084$ & $\lil 0.097363$ & $\lil{2}$ \\
 $\lil 20$ & $\lil 4814320320$ & $\lil 43021526040$ & $\lil 111304219860$ & $\lil 100872302880$ & $\lil 0.022091$ & $\lil{4/3}$ \\
 $\lil 22$ & $\lil 5454179550$ & $\lil 39892554000$ & $\lil 93242799000$ & $\lil 81714578400$ & $\lil 0.025027$ & $\lil{10/9}$ \\
 $\lil 24$ & $\lil 4073954144$ & $\lil 40186134868$ & $\lil 126323098182$ & $\lil 162790595856$ & $\lil 0.018694$ & $\lil{2}$ \\
 $\lil 26$ & $\lil 918069454$ & $\lil 12091107788$ & $\lil 51322797162$ & $\lil 88711954896$ & $\lil 0.004213$ & $\lil{12/11}$ \\
 $\lil 28$ & $\lil 857901000$ & $\lil 12427836600$ & $\lil 55357035900$ & $\lil 98053394600$ & $\lil 0.003937$ & $\lil{6/5}$ \\
 $\lil 30$ & $\lil 535673924$ & $\lil 10415825728$ & $\lil 65248580472$ & $\lil 171951637976$  & $\lil 0.002458$ & $\lil{8/3}$ \\
 $\lil 32$ & $\lil 58664256$ & $\lil 1599900552$ & $\lil 13444986588$ & $\lil 46806142904$ & $\lil 0.000269$ & $\lil{1}$ \\
 $\lil 34$ & $\lil 69404898$ & $\lil 1684816476$ & $\lil 13621926834$ & $\lil 47836532832$ & $\lil 0.000318$ & $\lil{16/15}$ \\ 
 $\lil 36$ & $\lil 46346428$ & $\lil 1439916356$ & $\lil 14571970374$ & $\lil 64004385832$ & $\lil 0.000213$ & $\lil{2}$ \\ 
 $\lil 38$ & $\lil 7381190$ & $\lil 318303280$ & $\lil 4219159800$ & $\lil 23451227440$ & $\lil 0.000034$ & $\lil{18/17}$ \\ 
 $\lil 40$ & $\lil 10176048$ & $\lil 359222796$ & $\lil 4396494114$ & $\lil 24594847992$ & $\lil 0.000047$ & $\lil{4/3}$ \\ 
 $\lil 42$ & $\lil 4153336$ & $\lil 201583172$ & $\lil 3188901438$ & $\lil 22696587504$ & $\lil 0.000019$ & $\lil{12/5}$ \\  
 $\lil 44$ & $\lil 526596$ & $\lil 37126032$ & $\lil 772483368$ & $\lil 6703381264$ & $\lil 0.000002$ & $\lil{10/9}$ \\ 
 $\lil 46$ & $\lil 291342$ & $\lil 21296376$ & $\lil 459181188$ & $\lil 4284667104$ & $\lil 0.000001$ & $\lil{22/21}$ \\ 
 $\lil 48$ & $\lil 239760$ & $\lil 19964064$ & $\lil 493227744$ & $\lil 5290003952$ & $\lil 0.000001$ & $\lil{2}$ \\ 
 $\lil 50$ & $\lil 91392$ & $\lil 7454520$ & $\lil 183370572$ & $\lil 2026286376$ & $\lil 4.2E-7$ & $\lil{4/3}$ \\ 
 $\lil 52$ & $\lil 8912$ & $\lil 1337188$ & $\lil 52081950$ & $\lil 819360400$ & $\lil 4.1E-8$ & $\lil{12/11}$ \\
 $\lil 54$ & $\lil 25320$ & $\lil 2992860$ & $\lil 97569690$ & $\lil 1348117880$ & $\lil 1.2E-7$ & $\lil{2}$ \\
 $\lil 56$ & $\lil 2952$ & $\lil 422196$ & $\lil 18140238$ & $\lil 326084664$ & $\lil 1.4E-8$ & $\lil{6/5}$ \\
 $\lil 58$ & $\lil 1654$ & $\lil 307068$ & $\lil 14158938$ & $\lil 264266960$ & $\lil 7.6E-9$ & $\lil{28/27}$ \\
 $\lil 60$ & $\lil 452$ & $\lil 110300$ & $\lil 6862242$ & $\lil 173593136$ & $\lil 2.1E-9$ & $\lil{8/3}$ \\
 $\lil 62$ & $\lil 26$ & $\lil 8248$ & $\lil 645804$ & $\lil 19784976$ & $\lil 1.2E-10$ & $\lil{30/29}$ \\
 $\lil 64$ & $\lil 48$ & $\lil 12528$ & $\lil 890688$ & $\lil 25971336$ & $\lil 2.2E-10$ & $\lil{1}$ \\
 $\lil 66$ & $\lil 24$ & $\lil 6744$ & $\lil 545796$ & $\lil 18824896$ & $\lil 1.1E-10$ & $\lil{20/9}$ \\ \hline
 \end{tabular}
\caption{ \label{G37Table} For the gaps that actually occur in $\pgap(\pml{37})$, this table lists 
the number of gaps and driving terms of length $j \le 4$.  
The gaps $g \ge 16$ have longer driving terms as well; the gap
 $g=66$ has driving terms up to length $16$.  Also tabulated are the current ratio $w_{g,1}(\pml{37})$ and the
 asymptotic value for this ratio. }
 \end{center}
 \end{table}
\renewcommand\arraystretch{1}

\subsection{Tracking the relative growth of the classes of gaps.}
The raw populations of gaps within the cycles of gaps $\pgap(\pml{p})$ grow by factors
of $(p-2)$.  This led us to look at the ratios $w_{g,j}(\pml{p}) = n_{g,j}(\pml{p}) / \prod (p-2)$.

Within each residue class modulo $10$ we will be adding up the ratios of an infinite number 
of gaps.  To compare these in a practical manner, we first note that due to the recursive construction
of $\pgap(\pml{p_k})$ from $\pgap(\pml{p_{k-1}})$, the closures within driving terms occur 
methodically.  As a result we see larger gaps introduced generally in order as the sieve progresses.
In Table~\ref{G37Table}, the column for $j=1$ illustrates this introduction of larger gaps.

As a first comparison of the distributions across the residue classes, we consider the average $\mu_h(w_g)$
as the size of the gaps within each class increases.
\begin{eqnarray*}
\left. \mu_0(w_g(\infty)) \right|_N & = & \frac{1}{N} \left[w_{10,1}(\infty) + w_{20,1}(\infty) + \cdots + w_{10N,1}(\infty) \right] \\
\left. \mu_2(w_g(\infty)) \right|_N & = & \frac{1}{N} \left[w_{2,1}(\infty) + w_{12,1}(\infty) + \cdots + w_{10N-8,1}(\infty) \right] \\
\left. \mu_4(w_g(\infty)) \right|_N & = & \frac{1}{N} \left[w_{4,1}(\infty) + w_{14,1}(\infty) + \cdots + w_{10N-6,1}(\infty) \right] \\
\left. \mu_6(w_g(\infty)) \right|_N & = & \frac{1}{N} \left[w_{6,1}(\infty) + w_{16,1}(\infty) + \cdots + w_{10N-4,1}(\infty) \right] \\
\left. \mu_8(w_g(\infty)) \right|_N & = & \frac{1}{N} \left[w_{8,1}(\infty) + w_{18,1}(\infty) + \cdots + w_{10N-2,1}(\infty) \right] 
\end{eqnarray*}
In Table~\ref{B10Table} we list the gaps $g < 100$ in their respective residue classes, along with each gap's 
asymptotic ratio and the average ratio for the class up to this gap.

\renewcommand\arraystretch{1.05}
\begin{table}
\begin{center}
\begin{tabular}{|ccr|ccr|ccr|} \hline
\multicolumn{9}{|c|}{gaps, asymptotic ratios, \& mean asymptotic ratios} \\ 
$g$ & $\lil w_g(\infty)$ & $\lil \mu_0( w_g(\infty))$ & $g$ & $\lil w_g(\infty)$ & $\lil \mu_2 (w_g(\infty))$ &
 $g$ & $\lil w_g(\infty)$ & $\lil \mu_4 ( w_g(\infty))$ \\ \hline
 & & $\lil 0.000$ & $\lil 2$ & $\lil 1$ & $\lil 1.000$ & $\lil 4$ & $\lil 1$ & $\lil 1.000$ \\
$\lil 10$ & $\lil 4/ 3$ & $\lil 1.333$ & $\lil 12$ & $\lil 2$ & $\lil 1.500$ & $\lil 14$ & $\lil 6 / 5$ & $\lil 1.100$ \\
$\lil 20$ & $\lil 4 /3$ & $\lil 1.333$ & $\lil 22$ & $\lil 10/ 9$ & $\lil 1.370$ & $\lil 24$ & $\lil 2$ & $\lil 1.400$ \\
$\lil 30$ & $\lil 8 / 3$ & $\lil 1.777$ & $\lil 32$ & $\lil 1$ & $\lil 1.277$ & $\lil 34$ & $\lil 16 / 15$ & $\lil 1.316$ \\ 
$\lil 40$ & $\lil 4 /3 $ & $\lil 1.666$ & $\lil 42$ & $\lil 12/5$ & $\lil 1.502$ & $\lil 44$ & $\lil 10 / 9$ & $\lil 1.275$ \\
$\lil 50$ & $\lil 4 /3$ & $\lil 1.600$ & $\lil 52$ & $\lil 12/ 11$ & $\lil 1.433$ & $\lil 54$ & $\lil 2$ & $\lil 1.396$ \\
$\lil 60$ & $\lil 8 / 3$ & $\lil 1.777$ & $\lil 62$ & $\lil 30/29$ & $\lil 1.376$ & $\lil 64$ & $\lil 1$ & $\lil 1.339$ \\
$\lil 70$ & $\lil 8 / 5$ & $\lil 1.752$ & $\lil 72$ & $\lil 2$ & $\lil 1.454$ & $\lil 74$ & $\lil 36 / 35$ & $\lil 1.300$ \\
$\lil 80$ & $\lil 4/3$ & $\lil 1.700$ & $\lil 82$ & $\lil 40/39$ & $\lil 1.406$ & $\lil 84$ & $\lil 12 / 5$ & $\lil 1.422$ \\
$\lil 90$ & $\lil 8 /3$ & $\lil 1.807$ & $\lil 92$ & $\lil 22/21$ & $\lil 1.370$ & $\lil 94$ & $\lil 46/45$ & $\lil 1.382$ \\ \hline \hline
$g$ & $\lil w_g(\infty)$ & $\lil \mu_6 ( w_g(\infty))$ & $g$ & $\lil w_g(\infty)$ & $\lil \mu_8 ( w_g(\infty))$ & & & \\ \hline
$\lil 6$ & $\lil 2$ & $\lil 2.000$ & $\lil 8$ & $\lil 1$ & $\lil 1.000$ & & & \\
$\lil 16$ & $\lil 1$ & $\lil 1.500$ & $\lil 18$ & $\lil 2$ & $\lil 1.500$ & & & \\
$\lil 26$ & $\lil 12/11$ & $\lil 1.363$ & $\lil 28$ & $\lil 6/5$ & $\lil 1.400$ & & & \\
$\lil 36$ & $\lil 2$ & $\lil 1.522$ & $\lil 38$ & $\lil 18/17$ & $\lil 1.314$ & & & \\
$\lil 46$ & $\lil 22/21$ & $\lil 1.427$ & $\lil 48$ & $\lil 2$ & $\lil 1.451$ & & & \\
$\lil 56$ & $\lil 6/5$ & $\lil 1.389$ & $\lil 58$ & $\lil 28/27$ & $\lil 1.382$ & & & \\
$\lil 66$ & $\lil 20/9$ & $\lil 1.508$ & $\lil 68$ & $\lil 16/15$ & $\lil 1.337$ & & & \\
$\lil 76$ & $\lil 18/17$ & $\lil 1.452$ & $\lil 78$ & $\lil 24/11$ & $\lil 1.443$ & & & \\
$\lil 86$ & $\lil 42/41$ & $\lil 1.404$ & $\lil 88$ & $\lil 10/9$ & $\lil 1.406$ & & & \\
$\lil 96$ & $\lil 2$ & $\lil 1.464$ & $\lil 98$ & $\lil 6/5$ & $\lil 1.385$ & & & \\ \hline
\end{tabular}
\caption{ \label{B10Table} The distribution of gaps $ g < 100$ that maintain pairs of last digits modulo $10$.}
 \end{center}
 \end{table}
\renewcommand\arraystretch{1}

Table~\ref{B10Table} helps us make a few observations about the effect of Theorem~\ref{PolThm} on the
average asymptotic ratios.  Gaps that are divisible by $3$ have a factor of $2$ in their asymptotic ratio, and
gaps that are divisible by $5$ have a factor of $4/3$ in theirs.  These are as large as these factors get.  If a gap
is divisible by a larger prime $p$, the asymptotic ratio $w_g(\infty)$ has a factor of $(p-1)/(p-2)$.

Working in base $10$, the gaps divisible by $3$ rotate through the classes in this order:  $h = 6, 2, 8, 4, 0$.
When the number $N$ of gaps in each class is small, we can see the impact of a mod-$3$ gap on 
the average.  In Table~\ref{B10Table}, look at $\mu_4$ jump at the gaps $g=24$ and $g=54$.  To make a fair
comparison across classes, we should pick $N$ to include a complete rotation of $3$ through the classes; e.g.
stopping when the gap $g$ is a multiple of $30$.

We also observe that the class $g = 0 \bmod 10$ will have all of the gaps divisible by $5$, giving the average for this
class a consistent boost.  Interestingly, the corresponding factor for $w_g(\infty)$ is $4/3$, which is the factor needed to 
compensate for this class having four ordered pairs $(a,b)$ as compared to the other classes having only
three ordered pairs of last digits.

The class $g = 0 \bmod 10$ will also contain all of the primorial gaps $g = \pml{p}$ for $p \ge 5$ and all of their multiples.  
By Theorem~\ref{PolThm} the primorial gaps mark new maxima for $w_g(\infty)$.  For the gap $g=30$, 
the asymptotic ratio is $w_{30}(\infty) = 8/3$, and for $g=210$ this asymptotic ratio jumps to
$w_{210}(\infty) = 48/15$.

Based on these asymptotic distributions, we expect that the biases observed by Oliver and Soundararajan will
disappear among large primes.  For the cycles of gaps $\pgap(\pml{p})$, the evolution of the dynamic system
plays out on massive scales.  Compared to the scales on which the primes evolve, the computations by Oliver and
Soundararajan \cite {OS} are very early.  For example, for their computations for distributions in base $3$ 
they considered the first million primes; these occur in the first twelfth of $\pgap(\pml{23})$ (that is, 
within the first two copies of $\pgap(\pml{19})$.
Their computations for distributions in base $10$ use the first one hundred million
primes; these occur in the first third of $\pgap(\pml{29})$.   We observe in Table~\ref{G37Table}
that even for the small gaps the ratios $w_{g,1}(\pml{37})$ at this stage are far from their asymptotic values.


\section{Rate of convergence to the ultimate distributions}\label{AsymSection}
From our observations above about the asymptotic ratios for small gaps, we believe that Oliver and
Soundararajan are observing transient phenomena.  In this section we demonstrate that these transient
biases will persist for any computationally tractable primes.

The asymptotic values described in the previous section evolve {\em very} slowly, even for small gaps.
For example, the gap $g=30$ ultimately occurs $4/3$ times as often as the gap $g=6$ and $8/3$ times as often
as the gap $g=2$; but the gap $g=30$ is not even more numerous than the gap $g=2$ until $\pgap(\pml{q})$
with $q \approx 2E6$.


From Equation~(\ref{EqEigSys}) we can determine more specifically the rate at which
the ratio $w_{g,1}(\pml{p})$ converges to the asymptotic value $w_{g,1}(\infty)$.
For computational simplicity, we use the approximation $a_j^k \approx (a_2^k)^{j-1}$
to express the expansion as a polynomial in $\lambda = a_2^k$.
\begin{equation*} 
w_{g,1}(\pml{p_k}) \approx L_1 \cdot w_g(\pml{p_0}) - (L_2 \cdot w_g(\pml{p_0})) \lambda +
 (L_3 \cdot w_g(\pml{p_0})) \lambda^2 - (L_4 \cdot w_g(\pml{p_0})) \lambda^3 + \cdots
\end{equation*}

This system converges to the asymptotic value $w_{g,1}(\infty)= L_1 w_g(\pml{p_0})$ as quickly as
$\lambda = a_2^k \fto 0$.  This decay of $a_2^k$ is extraordinarily slow.  The bounds in Equation~\ref{EqBnds}
show that this decay is proportionate to $1/ \ln p_k$.

To understand how far we are from convergence for a prime $q \approx 3E15$, let's again consider the gap $g=30$.
Ultimately the gap $g=30$ will have ratio $w_{30}(\infty) = 8/3$ compared to the gap $g=2$.
For $q \approx 3E15$, in $\pgap(\pml{q})$ the ratio of gaps $g=30$ to $g=2$ is $w_{30,1}(\pml{q}) \approx 1.976$.

For the next primorial gap $g=210$, the asymptotic ratio is $w_{210}(\infty) = 48/15$, but in $\pgap(\pml{q})$ we
only have $w_{210,1}(\pml{q}) \approx 0.265$.

For setting the vectors of initial values $w_g(\pml{p_0})$ we have an incentive to pick $p_0$ as large as possible, since
the dynamic system only holds exactly for those gaps $g$ all of whose prime factors are less than or equal to $p_0$.  Gaps
$g$ with larger prime factors will be underrepresented by factors of $(p-1)/(p-2)$; this is an application
of Corollary~5.7 in \cite{HRSFU}.

{\em Example.} To apply this polynomial approximation to the distribution of last digits of 
consecutive primes, we use $\pgap(\pml{37})$ to obtain initial conditions for the gaps $g=2,\ldots,420$.  This will give
accurate representations for all of the gaps that are not divisible by the larger primes $p = 41, 43, \ldots, 199$. 

\begin{figure}[t]
\centering
\includegraphics[width=5in]{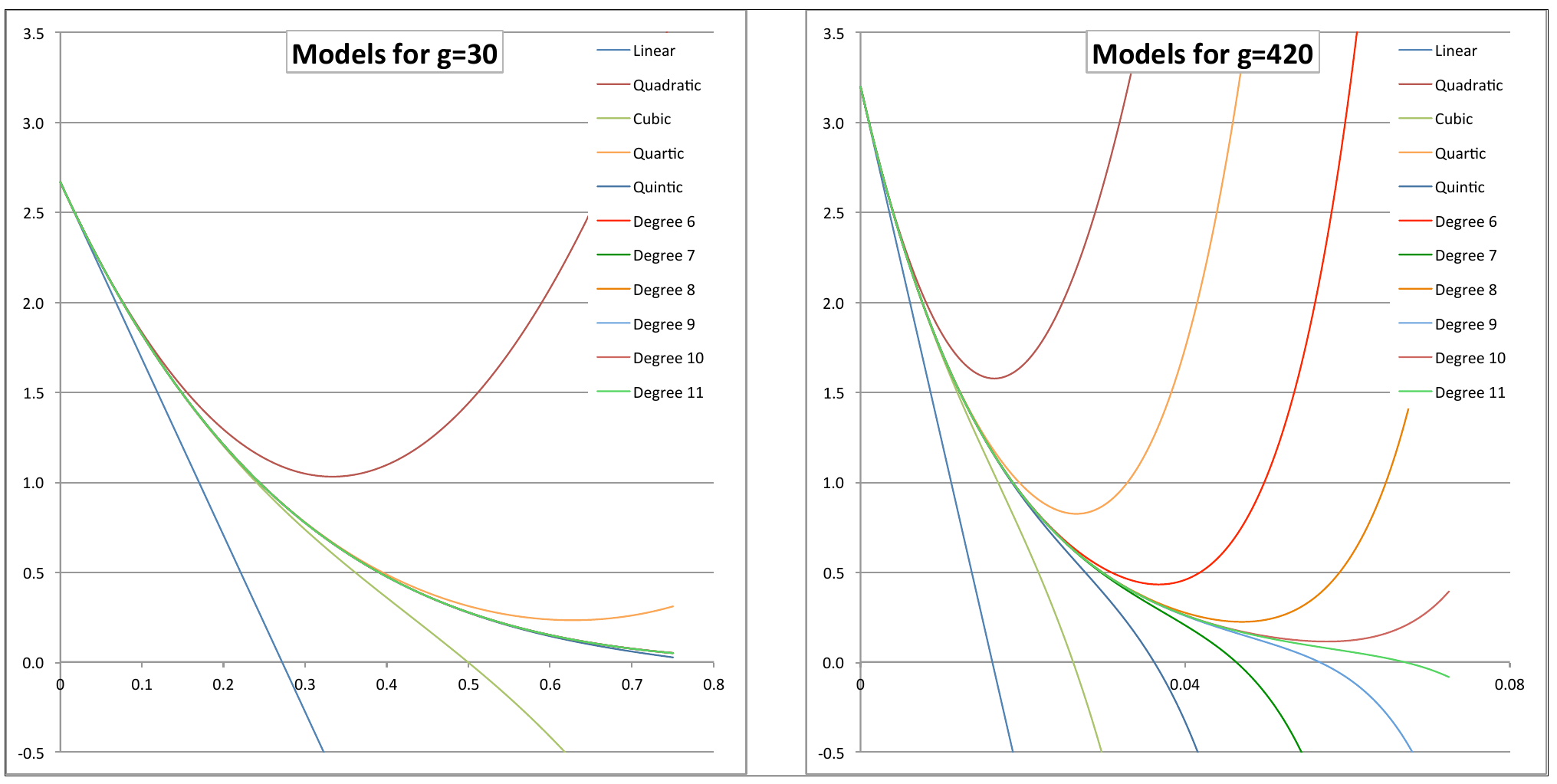}
\caption{\label{LPolyFig} Two examples of the polynomial approximations in Equation~\ref{EqEigSys}.  The approximations
differ from the exact discrete model by substituting $\lambda = (a_2^k)^j$ for $a_{j+1}^k$.  The gap $g=30$ has driving
terms up to length $j=8$, so the approximations of degree $8$ and higher coincide with that of degree $7$. }
\end{figure}

We take the first twelve terms of the polynomial approximation and apply 
this model to the ratios $w_{g,1}(\pml{p_k})$ for the gaps
$g=2,\ldots,420$ sorted into their residue classes modulo $10$.
For each gap $g$, we define the coefficients $l_i = L_i \cdot w_g(\pml{37})$.
Then we have the degree-11 model
\begin{equation*}
w_{g,1}(\pml{p_k})  \approx  l_1 - l_2 \lambda + l_3 \lambda^2 - \cdots
 + l_{11} \lambda^{10} - l_{12} \lambda^{11}
\end{equation*}
Graphs of these models for the gaps $g=30$ and $g=420$ from degree 1 to degree 11 are plotted
in Figure~\ref{LPolyFig}.  These graphs provide a sense of the range of values of $\lambda = a_2^k$
over which we can rely on the accuracy of these approximations.

Turning our attention back to the residue classes $\bmod 10$, 
we consider the aggregate model for the gaps in each residue class.
These five aggregate models are depicted in Figure~\ref{WratFig}.
As $p_k \fto \infty$, the parameter 
\[ \lambda = a_2^k = \prod_{p_1}^{p_k} \frac{p-3}{p-2} \fto 0. \]
Figure~\ref{WratFig} depicts the evolution of the populations
for each residue class versus $\log p_k$.  For each value of $p_k$, we normalize the total population of
gaps in class $h \bmod 10$ by the population of gaps with residue $2 \bmod 10$:
\[ W_h (\pml{p_k}) = \frac{\sum_{g = h \bmod 10} w_{g,1}(\pml{p_k})} 
{\sum_{g = 2 \bmod 10} w_{g,1}(\pml{p_k})}  \]

To give a sense of the initial conditions underlying this slow evolution of the populations of gap,
in $\pgap(\pml{37})$ the gap $g=420$ has a total of 
$697373938800$ driving terms.  These range from $2$ driving terms of length $j=47$ through
$304$ driving terms of length $j=76$.  Under the twelve-term approximation, the gap $g=420$ will not
reach $10\%$ of its asymptotic ratio of $w_{420,1}(\infty)=3.2$ until $a_2^k < 0.0365$; that is until
$\pgap(\pml{p_k})$ with $p_k \approx 1.12E45$.  With this model, we don't expect the gap $g=420$ to
be more numerous than the gap $g=2$ until $a_2^k < 0.01415$, that is until $p_k \approx  3.57E87$.

In Figure~\ref{WratFig} we can see how the biases observed by Oliver and Soundararajan will be corrected
for large primes.  The distributions calculated by Oliver and Soundararajan \cite{OS} correspond to the dashed
line toward the 
left of Figure~\ref{WratFig}.  Summing up their computed values by residue class we calculate the following
ratios for the first $10^8$ prime numbers:
 
\begin{center}
\begin{tabular}{|c|c|rr|} \hline
\multicolumn{4}{|c|}{Ratios $W_h$ from first $10^8$ primes \cite{OS} } \\ \hline
$h$ & $(a,b)$ & $\sum_g n_{g,1}$ & $W_h$  \\ \hline
$2$ & $(1,3), \; (7,9), \; (9,1)$ & $22852739$ & $1$  \\
$4$ & $(3,7), \; (7,1), \; (9,3)$ & $19790617$ & $0.866006$  \\
$6$ & $(1,7), \; (3,9), \; (7,3)$ & $21762703$ & $0.952302$  \\
$8$ & $(1,9), \; (3,1), \; (9,7)$ & $17466066$ & $0.764288$  \\
$0$ & $(1,1), \; (3,3), \; (7,7), \; (9,9)$ & $18127875$ & $0.793247 $ \\ \hline
\end{tabular}
\end{center}
The calculated ratios $W_h$ for the first $10^8$ primes are consistent with the distributions of gaps in 
$\pgap(\pml{p})$ for very small primes $p$.  Intuitively, we see that for small primes $p_k$ the residue classes
$h=2,4,6$ all start with significant populations of the gaps $g=2,4,6$ respectively, while the classes $h=8,0$ have
to manufacture representative gaps as the recursion on the cycle of gaps $\pgap(\pml{p_k})$ proceeds.

We see in Figure~\ref{WratFig} that these biases will change for much larger
primes.  Ultimately, for the sample of gaps $2 \le g \le 420$, these ratios will converge to:
\[ \begin{array}{clccl}
W_2(\infty) & =  1 & & W_8(\infty) & = 1.0026 \\
W_4(\infty) & = 1.0007 & & W_0(\infty) & = 1.3192 \\
W_6(\infty) & = 1.0029 & & & \\
\end{array}
\]

For this sample of gaps, the residue class $h = 0 \bmod 10$ has about $98.9\%$ of the size necessary to
compensate for having four corresponding ordered pairs $(a,b)$ as compared to the three for each
of the other residue classes.


\section{Distributions in other bases}\label{SecBases}
The work above has addressed the residue classes of primes in base $10$.  For base $10$ the dynamic system that
models populations of gaps across stages of Eratosthenes sieve indicates that the biases calculated by Oliver and 
Soundararajan are transient phenomena.  These biases will gradually fade for very large primes.  How is this analysis affected
by the choice of base?

Our work in base $10$ consisted of three components:  identifying the ordered pairs of last digits $(a,b)$ and gaps $g$
that correspond to each residue class; comparing the asymptotic ratios $w_{g,1}(\infty)$ for the gaps within each 
residue class; and looking at the initial conditions and rates of convergence for the gaps within each residue class.
As we consider other bases, we look at the effect that a new basis has on each of these components.

\begin{figure}[t]
\centering
\includegraphics[width=5in]{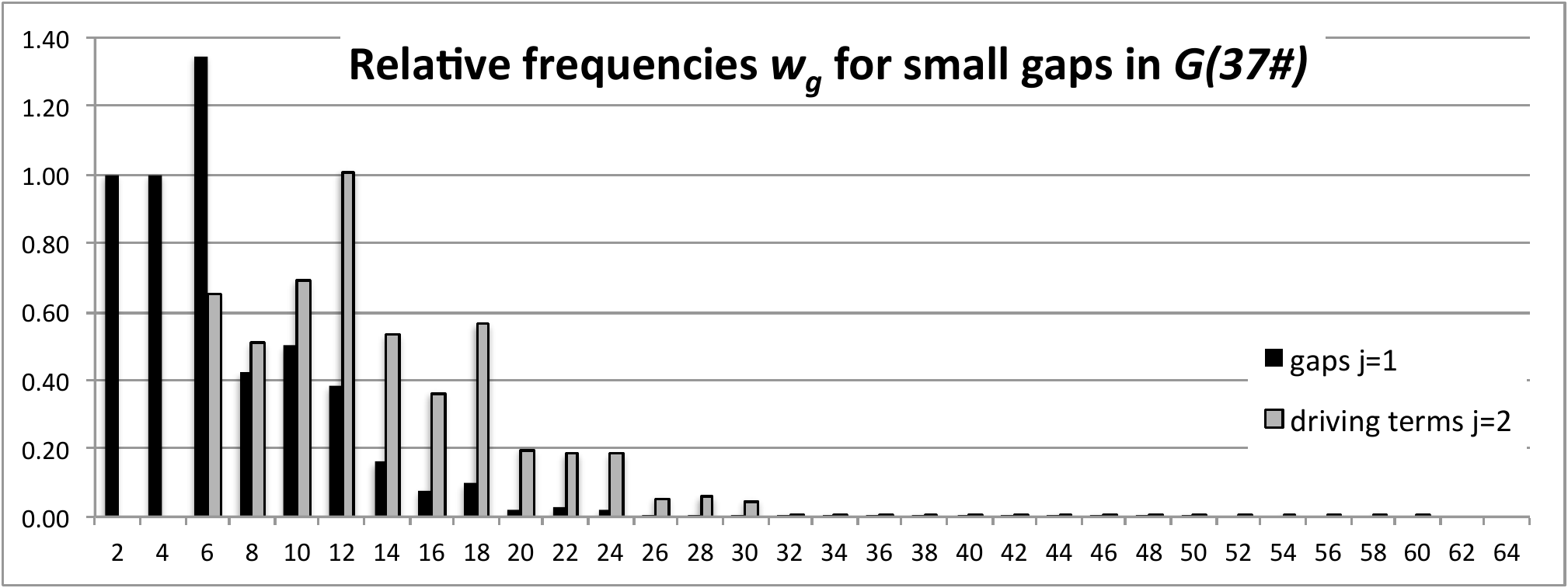}
\caption{\label{wG37Fig} The initial conditions for the gaps that occur in $\pgap(\pml{37})$. 
The relative frequencies 
$w_{g,j}(\pml{37})$ are shown for the gaps and for their driving terms of length $2$.
We see the early leads that the gaps $2,4,6$ enjoy.  Other small gaps, especially those divisible by $6$ have
short driving terms that will quickly boost their populations as well.}
\end{figure}

Since the asymptotic ratios tend toward uniformity across the ordered pairs, the choice of base will simply shift
the distributions around the residue classes.  The asymptotic ratios will be proportional to the number of ordered
pairs $(a,b)$ corresponding to that class.

By setting a base, we set the assignment of gaps, especially the small gaps, to the respective residue classes.
These residue classes inherit the initial biases and rates of convergence associated with the assigned gaps.

Note that the bias is dominated by small gaps, especially the gaps $2$, $4$ and $6$.  
Figure~\ref{wG37Fig} illustrates the components of the bias
introduced by small gaps.  We can see that the gaps
$g=12,10,18$ will quickly make significant contributions as well, and that for powers of $2$, e.g. $g=8,16,32$,
the populations will lag
compared to gaps of similar size.
The biases in the initial populations for small gaps, for example as illustrated in Table~\ref{G37Table} for $\pgap(\pml{37})$, 
will be inherited by the residue classes to which the small gaps belong.  

For the bases $3$ and $6$, all gaps that are multiples of $6$ fall into a single residue class.  
Similarly for the bases $5$ and $10$, all gaps that are multiples of $10$ fall into a single residue class.
The base $30$ separates the multiples of $6$ and $10$ into a small set of residue classes.  In contrast,
powers of $2$, like the base $8$, distribute the multiples of $3$ and $5$ (and any odd prime) across residue
classes. 

We illustrate this assignment of the initial bias with the gaps $2 \le g \le 420$ under the bases $3$, $8$, and $30$.

\subsection{Distributions in bases $3$ or $6$.}  Oliver and Soundararajan calculated the distributions of select pairs
$(a,b)$ modulo $3$ up through $10^{12}$, and they compare these favorably to a conjectured model derived
from the Hardy and Littlewood's \cite{HL} work on the $k$-tuple conjecture.

In our approach through $\pgap(\pml{p})$, for each residue class $h \bmod 3$ we identify the associated gaps $g$
and the ordered pairs with $h = b-a$ \cite{OS}.  We then calculate both the initial estimate $W_h(\pml{1993})$
based on the sample of gaps $g = 2, \ldots, 420$ and initial conditions from $\pgap(\pml{37})$; and the
asymptotic ratio $W_h(\infty)$ for this sample of gaps.  For the ratios, we normalize by the values for the
class $h = 2 \bmod 3$.

\begin{tabular}{|cl|c|rr|} \hline
$h \bmod {\bf 3}$ & $g$'s & $(a,b)$ & $W_h(\pml{1993})$ & $W_h(\infty)$ \\ \hline
$2$ & $2, 8, 14, \ldots$ & $(2,1)$ &  $1$ &  $1$ \\
$1$ & $4, 10, 16, \ldots$ & $(1,2)$ &  $1.0009$ &  $1.0010$ \\
$0$ & $6, 12, 18, \ldots$ & $(1,1), (2,2)$ & $1.6358$ & $1.9868$ \\ \hline
\end{tabular}

We see that the ratios $W_1(\pml{1993})$ and $W_0(\pml{1993})$ are consistent with Oliver and Soundararajan's
tabulations for base $3$.  Additionally, the asymptotic ratios for the sample of gaps indicates that the initial
bias will again disappear.

We note that these results for the primes modulo $3$ (or in base $3$) can be translated directly
into results for the primes modulo $6$ (or in base $6$).  Indeed, for any odd base $B$, there is a direct
translation of the residue classes, gaps, and ordered pairs into the base $2B$.  For example, Oliver and Soundararajan's 
computations for base $5$ could be combined with their initial computations base $10$.

One interesting aspect of working in base $3$ or $6$ is that all multiples of $6$ will fall in the class $h = 0 \bmod 6$,
and thus all of the primorials $g=\pml{p}$ will fall within this class. 
In $\pgap(\pml{11})$ the gap $g=6$ is the most frequent gap, and it grows more quickly than other gaps for many more
stages of the sieve. 

\subsection{Distributions in base $8$.}  
In base $8$ the small gaps are distributed more evenly across the residue classes.  We observe that 
the residue class $h = 0 \bmod 8$ starts slowly, and for this sample of gaps $2 \le g \le 420$ this residue
class still lags in its asymptotic value.  

\begin{tabular}{|cl|c|rr|} \hline
$h \bmod {\bf 8}$ & $g$'s & $(a,b)$ & $W_h(\pml{1993})$ & $W_h(\infty)$ \\ \hline
$2$ & $2, 10, 18, \ldots$ & $\lil (1,3), (3,5), (5,7), (7,1),$ &  $1$ &  $1$ \\
$4$ & $4, 12, 20, \ldots$ & $\lil (1,5), (5,1), (3,7), (7,3)$ &  $0.9695$ &  $1.0185$ \\
$6$ & $6, 14, 22, \ldots$ & $\lil (1,7), (7, 5), (5,3), (3, 1)$ &  $1.0086$ &  $1.0003$ \\
$0$ & $8, 16, 24, \ldots$ & $\lil (1,1), (3,3), (5,5), (7,7)$ &  $0.7081$ &  $0.9676$ \\ \hline
\end{tabular}

\subsection{Distributions in base $30$.}
The next primorial base is $30 = \pml{5}$.  This base is big enough that the small gaps are well separated, and the multiples
of $g=6$ and $g=10$ fall into a few distinct classes.  The early bias toward small gaps $g=2,4,6,10,12$ and 
even $g=14,18$ fall into separate residue classes.  We see the early biases in $W_h(\pml{1993})$ 
for base $30$ in Table~\ref{B30Table}. 

\begin{table}
\begin{tabular}{|cl|c|rr|} \hline
$h \bmod {\bf 30}$ & $g$'s & $(a,b)$ & $W_h(\pml{1993})$ & $W_h(\infty)$ \\ \hline
$\lil 2$ & $\lil 2, 32, \ldots$ & $\lil (29,1), (11,13), (17,19)$ &  $\lil 1$ &  $\lil 1$ \\
$\lil 4$ & $\lil 4, 34, \ldots$ & $\lil (7,11), (13,17), (19,23)$ &  $\lil 1.0180$ &  $\lil 1.0019$ \\
$\lil 6$ & $\lil 6, 36, \ldots$ & $\lil (1,7), (7,13), (13,19)$ &  $\lil 1.7771$ &  $\lil 2.0021$ \\
 & & $\lil (11,17), (17,23), (23,29)$ & & \\
$\lil 8$ & $\lil 8, 38, \ldots$ & $\lil (11,19), (23,1), (29,7)$ &  $\lil 0.8154$ &  $\lil 1.0000$ \\
$\lil 10$ & $\lil 10, 40, \ldots$ & $\lil (1,11), (7,17), (13,23), (19,29)$ &  $\lil 1.0421$ &  $\lil 1.3245$ \\
$\lil 12$ & $\lil 12, 42, \ldots$ & $\lil (1,13), (7,19), (11,23),$ &  $\lil 1.4228$ &  $\lil 1.9918$ \\
 & & $\lil (17,29), (19,1), (29,11) $ & & \\
$\lil 14$ & $\lil 14, 44, \ldots$ & $\lil (17,1), (23,7), (29,13)$ &  $\lil 0.7501$ &  $\lil 1.0028$ \\
$\lil 16$ & $\lil 16, 46, \ldots$ & $\lil (1,17), (7,23), (13,29)$ &  $\lil 0.5890$ &  $\lil 1.0015$ \\
$\lil 18$ & $\lil 18, 48, \ldots$ & $\lil (1,19), (11,29), (13,1),$ &  $\lil 1.0775$ &  $\lil 1.9956$ \\
 & & $\lil (19,7), (23,11), (29,17)$ & & \\
$\lil 20$ & $\lil 20, 50, \ldots$ & $\lil (11,1), (17,7), (23,13), (29,19)$ &  $\lil 0.6116$ &  $\lil 1.3287$ \\
$\lil 22$ & $\lil 22, 52, \ldots$ & $\lil (1,23), (7,29), (19,11)$ &  $\lil 0.5109$ &  $\lil 1.0020$ \\
$\lil 24$ & $\lil 24, 54, \ldots$ & $\lil (7,1), (13,7), (19,13)$ &  $\lil 0.8031$ &  $\lil 1.9920$ \\
 & & $\lil (17,11), (23,17), (29,23)$ & & \\
$\lil 26$ & $\lil 26, 56, \ldots$ & $\lil (11,7), (17,13), (23,19)$ &  $\lil 0.3920$ &  $\lil 1.0019$ \\
$\lil 28$ & $\lil 28, 58, \ldots$ & $\lil (1,29), (13,11), (19,17)$ &  $\lil 0.4122$ &  $\lil 1.0086$ \\
$\lil 0$ & $\lil 30, 60, \ldots$ & $\lil (1,1), (7,7), (11,11), (13,13),$ & $\lil 0.7578$ & $\lil 2.6153$ \\ 
 & & $\lil (17,17), (19,19), (23,23), (29,29)$ & & \\ \hline
\end{tabular}
\caption{\label{B30Table} The table for the distribution in base $30$.}
\end{table}


\section{Conclusion}
By identifying structure among the gaps in each stage of Eratosthenes sieve, we 
have been able to develop an exact model for the populations of gaps and their
driving terms across stages of the sieve.  
We have identified a model for a discrete dynamic system that takes 
the initial populations of a gap $g$ and all its driving terms in a cycle of gaps 
$\pgap(\pml{p_0})$ such that $g < 2p_1$, and thereafter provides the exact
populations of this gap and its driving terms through all subsequent cycles of gaps.

All of the gaps between primes are generated out of these cycles of gaps, with the gaps at the front of 
the cycle surviving subsequent closures.  The trends across the stages of Eratosthenes sieve indicate
probable trends for gaps between primes.  We are not yet able to translate the precision of the model for
populations of gaps in $\pgap(\pml{p})$ into a robust analogue for gaps between primes.

For the first $10^8$ primes, Oliver and Soundararajan \cite{OS,OSQ} calculated how often the possible pairs
$(a,b)$ of last digits of consecutive primes occurred, and they observed biases.  Regarding their calculations
they raised two questions: Does the observed bias persist?  Is the observed bias dependent upon the base? 
We have addressed  both of these questions by using the dynamic system that exactly models the populations
of gaps across stages of Eratosthenes sieve. 

The observed biases are transient phenomena.  The biases persist through the range of computationally tractable 
primes.  The asymptotics of the dynamic system play out on superhuman scales -- for example, continuing Eratosthenes
sieve at least through all $16$-digit primes.  To put this in perspective, the cycle $\pgap(\pml{199})$ has more gaps than there are
particles in the known universe; yet in $\pgap(\pml{p})$ for a $16$-digit prime $p$, small gaps like $g=30$ will still be appearing
in frequencies well below their ultimate ratios.  Gaps the size of $g=210$ will just be emerging, relative to the
prevailing populations of small gaps at this stage.

Our work on the relative frequency of gaps modulo $10$ for $\pgap(\pml{p})$ has addressed the bias between the
residue classes.  The observed biases are due to the quick appearance of small gaps and the slow evolution
of the dynamic system.
While we have addressed the inter-class bias, we have said nothing about the intra-class bias, that is, unequal distributions
across the ordered pairs $(a,b)$ within a given residue class modulo $10$.  Our initial calculations here indicate that
this bias should also disappear eventually, but this exploration needs to be more thorough.
The model developed by Oliver and Soundararajan also depends only on the residue class $h=b-a$.

Our calculations use a sample of gaps $g=2,\ldots,420$.  To improve the precision of our calculations of the asymptotic
ratios $w_h(\infty)$ across residue classes, it would be useful to find a normalization that makes working with all gaps
$g=2n$ manageable.

Once we understand the model for gaps, then any choice of base reassigns the gaps across the residue classes for this
base.  The number of ordered
pairs corresponding to a residue class $h$ is proportional to the asymptotic relative frequency $W_h(\infty)$.  
The initial biases and more
rapid convergence that favor the small gaps can be observed, over any computationally tractable range,
for the residue classes to which these small
gaps are assigned.

\bibliographystyle{amsplain}

\providecommand{\bysame}{\leavevmode\hbox to3em{\hrulefill}\thinspace}
\providecommand{\MR}{\relax\ifhmode\unskip\space\fi MR }
\providecommand{\MRhref}[2]{%
  \href{http://www.ams.org/mathscinet-getitem?mr=#1}{#2}
}
\providecommand{\href}[2]{#2}

\end{document}